# Axiomatic Theory of Algorithms:

# Computability and Decidability in Algorithmic Classes

**Mark Burgin**

Department of Mathematics
University of California, Los Angeles
405 Hilgard Ave.
Los Angeles, CA 90095

**Abstract**

Axiomatic approach has demonstrated its power in mathematics. The main goal of this paper is to show that axiomatic methods are also very efficient for computer science. It is possible to apply these methods to many problems in computer science. Here the main modes of computer functioning and program execution are described, formalized, and studied in an axiomatic context. The emphasis is on three principal modes: *computation, decision,* and *acceptation*. Now the prevalent mode for computers is *computation*. Problems of artificial intelligence involve *decision mode*, while communication functions of computer demand *accepting mode*. The main goal of this paper is to study properties of these modes and relations between them. These problems are closely related to such fundamental concepts of computer science and technology as computability, decidability, and acceptability. In other words, we are concerned with the question what computers and software systems can do working in this or that mode. Consequently, results of this paper allow one to achieve higher understanding of computations and in such a way, to find some basic properties of computers and their applications. Classes of algorithms, which model different kinds of computers and software, are compared with respect to their computing, accepting or deciding power. Operations with algorithms and machines are introduced. Examples show how to apply axiomatic results to different classes of algorithms and machines in order to enhance their performance.

# 1 Introduction

The world of computers and their applications is very complex and sophisticated. It involves interaction of many issues: social and individual, biological and psychological, technical and organizational, economical and political. However, humankind in its development created a system of intellectual "devices" for dealing with overcomplicated systems. This system is called science and its "devices" are theories.

When people want to see what they cannot see with their bare eyes, they build and use various magnifying devices. To visualize what is situated very far from them, people use telescopes. To discern very small things, such as microbes or cells of living organisms, people use microscopes. In a similar way, theories are "magnifying devices" for mind. They may be utilized both as microscopes and telescopes. Being very complex these "theoretical devices" have to be used by experts.

Complexity of the world of modern technology is reflected in a study of Gartner Group's TechRepublic unit (Silverman, 2000). According to it, about 40% of all internal IT projects are canceled or unsuccessful, meaning that an average of 10% of a company's IT department each year produces no valuable work. An average canceled project is terminated after 14 weeks, when 52% of the work has already been done, the study shows. In addition, companies spend an average of almost $1 million of their $4.3 million annual budgets on failed projects, the study says. However, companies might be able to minimize canceled projects as well as the time for necessary cancellation if they have relevant evaluation theory and consult people who know how to apply this theory.

All developed theories have a mathematical ground core. Thus, mathematics helps science and technology in many ways. Scientists are even curious, as wrote the Nobel Prize winner Eugene P. Wigner in 1959, why mathematics being so abstract and remote from reality is unreasonably effective in the natural sciences. It looks like a miracle.

So, it is not a surprise that mathematics has its theories for computers and computations. The main of these theories is *theory of algorithms*. It explains in a logical way how computers function and how they are organized. Here, we are going to investigate functioning of computers, their networks, and other information processing systems.

To do this, three main modes of computer functioning are separated: *computation, decision,* and *acceptation*. The prevalent mode now is *computation*. Computers were developed with the goal to compute. However, the development of computer technology and emergence of new kinds of information processing systems essentially extended functions of computers. As a result, solving different problems has become the main goal of computers. This involves mostly *decision mode* of functioning. Decision-making is a prime goal of artificial intelligence. At the same time, appearance of Internet transformed computers into communication devices. Communication involves receiving and sending information. Sending information is realized in *computing mode*, while receiving information demands *accepting mode*. For example, one of the vital problems for computer security is to make a decision whether to accept a message or to reject it because this message contains a virus. The main goal of this paper is to study properties of modes of computer functioning and relations between them.

Usually to study properties of computers and to develop more efficient applications, we use mathematical models. There is a variety of such models: Turing machines of different kinds (with one tape and one head, with several tapes, with several heads, with n-dimensional tapes, non-deterministic, probabilistic, and alternating Turing machines, Turing machines that take advice and Turing machines with oracle, etc.), Post productions, partial recursive functions, neural networks, finite automata of different kinds (automata without memory, autonomous automata, accepting automata, probabilistic automata, etc.), Minsky machines, normal Markov algorithms, Kolmogorov algorithms, formal grammars of different kinds (regular, context free, context sensitive, phrase-structure, etc.), Storage Modification Machines or simply, Shönhage machines, Random Access Machines (RAM), Petri nets, which like Turing machines have several forms – ordinary, regular, free, colored, self-modifying, etc.), and so on.

This diversity of models is natural and useful because each of these classes is suited for some kind of problems. In other words, the diversity of problems that are solved by computers involves a corresponding diversity of models. For example, general problems of computability involve such models as Turing machines and partial recursive functions. Finite automata are used for text search, lexical analysis, and construction of semantics for programming languages. In addition, different computing devices demand corresponding mathematical models. For example, universal Turing machines and inductive Turing machines allows one to investigate characteristics of conventional computers (Burgin, 2001). Petri nets are useful for modeling and analysis of computer networks, distributed computation, and communication processes (Peterson, 1981). Finite automata model computer arithmetic. Neural networks reflect properties of the brain. Abstract vector and array machines model vector and array computers (Pratt *et al*,1974)

To utilize some kind of models that are related to a specific type of problems, we need to know its properties. In many cases, different classes have the same or similar properties. As a rule, such properties are proved for each class separately. Thus, alike proofs are repeated many times in similar situations involving various models and classes of algorithms.

In contrast to this, the axiomatic theory of algorithms suggests a different approach. Assuming some simple basic conditions, we derive in this theory many profound properties of algorithms. This allows one, when dealing with a specific model not to prove this property, but only to check the conditions from the assumption, which is much easier than to prove the property under consideration. Thus, we can obtain various characteristics of types of computers and software systems.

In addition, the axiomatic theory of algorithms solves another problem related to the very concept of algorithm. Now there is no consent on the definition of algorithm. Mathematicians and computer scientists are still asking what algorithm is (cf., for example, (Moschovakis, 2001)). Different researches give their own answer to this question. The axiomatic theory of algorithms allows one to derive properties of algorithms without exact specification of the concept of algorithm. Thus, its results are also useful for different concepts considered in theory of algorithm. As the result, this paper develops further the axiomatic approach originated in (Burgin, 1985).

It is necessary to remark that logical tools and axiomatic description has been used in computer science for different purposes. For example, Manna (1974) constructed axiomatic theory of programs, while Milner (1989) developed axiomatic theory of communicating processes. All such approaches described by axioms separated objects. For instance, in the theory of Manna, such objects are computer programs, while in the theory of Milner, such objects are communicating computational processes. Consequently, they use local axiomatization. In comparison with those approaches, the approach presented in this paper is global as axioms are used to describe classes of algorithms, programs or automata.

## 2 Computational Modes in the Axiomatic Context

The axiomatic approach unifies studies of diverse models and classes of algorithms. Different models of algorithms represent, as a rule, distinct classes of programs and computational devices such as computers, calculators and transducers. Different classes of algorithms within the same model represent, in general, individual classes of programs and processes that are realized by the same or similar computational devices. For example, working with such model as Turing machines, we can fix the number of tapes and/or heads. This allows us to model hardware of a computer. Taking different systems of rules for these Turing machines, we pull out specific classes of Turing machines and model in such a way separate programs with which the computer is working. Making restrictions on time of the Turing machine functioning, we take into account specific classes of problems and investigate in such a way real-time or other tractable computations. Here we apply the axiomatic approach to classification and study of the modes of computer functioning. Such general setting allows one to consider properties of computational schemes and algorithms independently of types and kinds of computers and their software.

In addition, we show how these modes are related to general concepts of computability and decidability. To separate main modes in the axiomatic context, we assume validity of two basic axioms for algorithms (Burgin, 1985).

**A1.** Any algorithm $A$ determines a binary relation $r_A$ in the direct product $X \times Y$ of all its inputs $X$ and all its outputs $Y$.

**Remark 2.1.** Usually this relation $r_A$ is defined by computability in a corresponding class of algorithms or by computability by $A$. That is, a pair $(u, z)$ from $X \times Y$ belongs to $r_A$ if and only if application of $A$ to $u$ results in $z$. However, algorithms without output are also considered and utilized. Accepting finite and pushdown automata may be taken as examples. Usually it is assumed that such abstract machines do not give an output (see, for example, Hopcroft *et al*, 2001). Actually these automata give as an output messages "accepted" and "rejected." In such a way, it is possible to correspond to an algorithm $A$ those relations that are acceptable (cf. Definition 2.4) or decidable (cf. Definition 2.7) by $A$.

Axiom A1 may be formulated in a more restricted version for deterministic algorithms.

**A1d.** Any deterministic algorithm $A$ determines a partial function $f_A$ from the set $X$ of all its inputs to the set $Y$ of all its outputs.

**Remark 2.2.** Axiom A1d is a necessary but not a sufficient condition for an algorithm to be deterministic. For example, nondeterministic finite automata also define an output function of acceptance, which has two values "accepted" and "rejected."

Usually this function $f_A$ is defined by computability. That is, the value of $f_A$ at the point $u$ from the set $X$ is equal to the result $A(u)$ (if it exists) of applying $A$ to $u$. However, it is possible to correspond to $A$ those functions $f_A$ that are defined by acceptability through $A$ (cf. Definition 2.4) or by decidability by $A$ (cf. Definition 2.7).

**Remark 2.3.** Instead of introducing axiom A1d of determinism, we can simply define a deterministic algorithm $A$ as an algorithm for which the relation $r_A$ is a partial function.

**Definition 2.1.** Two algorithms are called functionally equivalent with respect to computability (acceptability, positive decidability, negative decidability or decidability) if they define in the corresponding mode the same function $f_A$ or relation $r_A$.

**Example 2.1.** In the theory of finite automata, functional equivalence means that two finite automata accept the same language (Hopcroft *et al*, 2001). This relation is used frequently to obtain different properties of finite automata. The same is true for the theory of pushdown automata.

Algorithms that work with finite words in some alphabet $X$ are the most popular in theory of algorithms. As a rule, only finite alphabets are utilized. For example, natural numbers in the decimal form are represented by words in the alphabet $X = \{0, 1, 2, 3, 4, 5, 6, 7, 8, 9\}$, while in binary form they are represented by words in the alphabet $X = \{0, 1\}$. The words in $X$ may represent natural numbers or not, but in any way there is a natural procedure to enumerate all such words. This makes it possible, when it is necessary, to assume that algorithms work with natural numbers. In such a way, through enumeration of words, any algorithm $A$ defines a partial function $f_A : \mathbf{N} \to \mathbf{N}$ (cf., (Burgin, 1985)). However, there are many reasons to consider algorithms that work with infinite words (Vardi and Wolper, 1994) or with such infinite objects as real numbers (Blum, *et al*. 1998)

**Remark 2.4.** Many algorithms (cf., for example, (Krinitsky, 1977) or (Burgin, 1985)) work with more general entities than words. As an example, we may consider configurations, which are utilized in (Kolmogorov, 1953) as inputs and outputs of algorithms. Configurations are sets of symbols connected by relations and may be treated as multidimensional words. Discrete graphs are examples of configurations.

A general idea of algorithm (cf., for example, (Rogers, 1987) or (Balcazar, Diaz, and Gabarro, 1988)) implies that there are three modes of processing input data:

1. *Computing* when an algorithm produces (computes or outputs) some words (its output or configuration) as a result of its activity.

2. *Deciding* when an algorithm, given a word/configuration $u$ and a set $X$ of words/configurations, indicates (decides) whether this word/configuration belongs to $X$.

3. *Accepting* when an algorithm, given a word/configuration $u$, either accepts this word/configuration or not.

These three types define not only the principal modes of computer functioning, but also the main utilization modes for algorithms and programs.

**Definition 2.2.** An algorithm $A$ accepts a word $u$ if $A$ gives a result when $u$ is its input.

**Example 2.2.** In the theory of finite automata, such acceptance is called acceptance by a result (Trahtenbrot and Barzdin, 1970). It is proved that acceptance by a result is functionally equivalent to acceptance by a state.

**Remark 2.5.** For many classes of algorithms or abstract automata, acceptance of a word $u$ means that the automata that works with the input $u$ comes to some inner state that is an accepting state for this algorithm or automata. Finite automata give an example of such a class. However, for such algorithms that produce some output, the acceptance assumption means that whenever an algorithm comes to an inner accepting state it produces some chosen result (e.g., the number 1) as its output. In such a way this algorithm informs that it has reached an inner accepting state.

Another way to define an accepting state is to consider a state of some component of an abstract automaton. For example, pushdown automata accept words not only by an accepting inner state, but also by an empty stack, that is, by a definite state of its stack.

Thus, given a class **A** of algorithms, we can separate specific types of sets that are determined by this class.

**Definition 2.3 (Computation).** A set $X$ is called computable in **A** if it is computable by some algorithm $A$ from **A**, i.e., if the output of $A$ (the range of $f_A$) is equal to $X$.

This is the usual way of computer functioning: from input to output.

**Definition 2.4 (Enumeration).** A set $X$ is called enumerable in **A** if it is computable by some algorithm $A$ from **A** and $A$ is defined on the whole $N$, i.e., $f_A$ is a total function.

Informally, enumeration means computation of all output values sequentially, one after another.

**Lemma 2.1.** Any set $X$ enumerable in **A** is computable in **A**.

Thus, we can see that enumeration is a particular case of computation.

**Example 2.3.** In the theory of Turing machines and recursive functions, the sets that are enumerable (computable) by Turing machines are called recursively enumerable (computable) sets. They play an important role. They constitute the largest class of sets that are constructible/computable by recursive algorithms. Application of the Church-Turing Thesis gives the conclusion that recursively enumerable sets constitute the largest class of sets that are algorithmically generated. This is based on the following conversion of Lemma 2.1.

**Proposition 2.1.** Any recursively computable set $X$ is recursively enumerable.

Proof. Here we give a semiformal proof of this statement. There are standard procedures (cf., for example, Ebbinghaus *et al*, 1970; or Rogers, 1987) that allow one to convert such proofs into exact mathematical deductions, in which all rules for Turing machines are given explicitly and each intermediate statement is deduced by a sequence of exact derivations.

Let us suppose that $X$ is a recursively computable infinite set. It means that there is a Turing machine $T$ that all outputs of $T$ constitute $X$. If $T$ is defined for all inputs, then $X$ is recursively enumerable and everything is proved. Otherwise, $T$ gives no output for some inputs, while for enumeration it is necessary to transform each input into some output. To remedy this deficiency and to preserve at the same time the initial output domain (the range of the function $f_T$), we build a Turing machine $D$ that realizes the process that is called the *bidiagonal covering* of $T$. Informally, the *bidiagonal covering* means that $D$ imitates all acts of computation that $T$ performs with all possible inputs. In particular, $D$ outputs all results of $T$ and only these results. This is proved by mathematical induction that is based on the schema of the Turing machine $D$, which is given in Figure 1. The description of functioning of $D$ is given in Table 1.

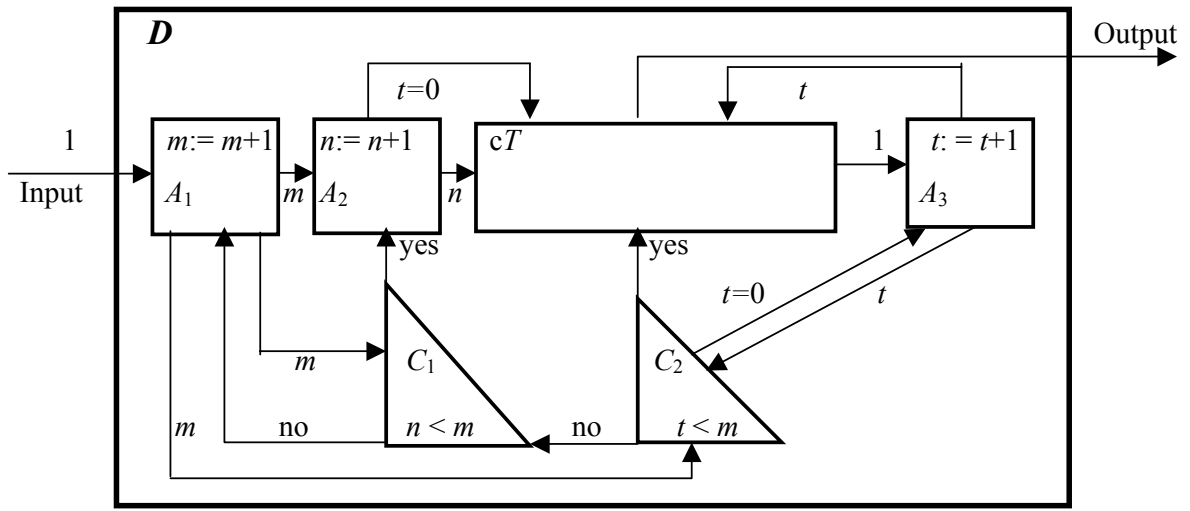

**Figure 1.** The structure of the Turing machine $D$ that realizes the diagonalization process for a Turing machine $T$; in it $cT$ is a copy of the machine $T$

1. $D$ initiates $A_1$ : $A_1$ takes 1 as its input and produces $m = 2$ as its output; then sends $n = 0$ as input to $A_2$ and $m = 2$ as input to $C_1$ and $C_2$ .

2. $D$ initiates $A_2$ : $A_2$ takes 0 as its input and produces $n = 1$ as its output; then sends $n = 1$ as input to $cT$ and $C_1$.

3. $D$ initiates $cT$ : $cT$ takes 1 as its input and simulates one step of $T$, obtaining as the result $T_1(1)$ as its output; if $cT$ comes to a final state of $T$ (but not of $D$), then $cT$ writes (on the output tape of $D$) the word $T_1(1)$ and the symbol * as the output of $D$; after this or if $cT$ does not come to a final state of $T$ , $cT$ sends 0 as input to $A_3$ .

4. $D$ initiates $A_3$ : $A_3$ takes 0 as its input and produces $t = 1$ as its output; then sends $t = 1$ as input to $C_2$ .

5. $D$ initiates $C_2$ : $C_2$ compares 1 and 2; as $1 < 2$, $C_2$ sends "yes" to $cT$ , opening $cT$ for one more step and goes to the stage 6.

6. $D$ initiates $cT$ : $cT$ simulates one more step of $T$, going from $T_t(n)$ to $T_{t+1}(n)$ as its output; if $cT$ comes to a final state of $T$ (but not of $D$), then $cT$ writes (on the output tape of $D$) the word $T_{t+1}(n)$ and the symbol * as the output of $D$; after this or if $cT$ does not come to a final state of $T$ , $cT$ sends 1 as input to $A_3$ .

7. $D$ initiates $A_3$ : $A_3$ takes 1 as its input and produces $t := t + 1$ as its output; then sends the new value of $t$ as input to $C_2$ .

8. $D$ initiates $C_2$ : $C_2$ compares $t$ and $m$; if $t < m$, then $C_2$ sends "yes" to $cT$ , opening $cT$ for one more step, and goes to the stage 6; if $t = m$, then $C_2$ sends $t = 0$ to $A_3$ , and goes to the stage 9.

9. $D$ initiates $C_1$ : $C_1$ compares $n$ and $m$; if $n < m$, then $C_2$ sends "yes" to $A_2$ , opening $A_2$ for one more step, and goes to the stage 11; if $n = m$, then $C_2$ sends "no" to $A_2$ , opening $A_2$ for one more step, and goes to the stage 10.

10. $D$ initiates $A_1$ : $A_1$ produces $m := m + 1$ as its output; then sends $n = 0$ as input to $A_2$ and the new value of $m$ as input to $C_1$ and $C_2$ .

11. $D$ initiates $A_2$ : $A_2$ produces $n := n + 1$ as its output; then sends and the new value of $n$ as input to $cT$ and $C_1$.

12. $D$ initiates $cT$ : $cT$ takes $n$ as its input and simulates one step of $T$, obtaining as the result $T_1(n)$ as its output; if $cT$ comes to a final state of $T$ (but not of $D$), then $cT$ writes (on the output tape of $D$) the word $T_1(n)$ and the symbol * as the output of $D$; after this or if $cT$ does not come to a final state of $T$ , $cT$ sends 0 as input to $A_3$ , and goes to the stage 7.

**Table 1:** The description of functioning of the Turing machine $D$ that realizes the diagonalization process for a Turing machine $T$

Using the *bidiagonal covering*, which is realized by the Turing machine $D$, we can correspond to $T$ a Turing machine $DT$ with the following property: if $T$ is defined for a number $m$, then $DT(n) = T(m)$ and if we take all smaller numbers 1, 2, ... , $m - 1$, then the Turing machine $T$ is defined exactly for $n - 1$ of them. This machine $DT$ works in the following manner. Given a number $n$ as an input, $DT$ uses $D$ to find $n$ first results of $T$. The machine $DT$ gives the last of these results as it output for the input $n$. By induction, we demonstrate that the Turing machine $DT$ is defined for any number $k < n$. These properties of the Turing machine $DT$ show that $DT$ defines a total function $f_{DT}$ and has the same range as $T$, By Definition 2.4, this means that $DT$ recursively enumerates the infinite set $X$.

Thus, we have proved that any recursively computable infinite set is recursively enumerable.

When the set $X$ is finite, we enhance the Turing machine $D_oT$ in the following way. After the first result $x$ of the machine $T$ is obtained by $D_oT$ in some cycle $n$, $D_oT$ gives $x$ as its result for inputs 1, 2, ... , $n$, i.e., $D_oT(i) = n$ for all $i = 1, 2, ... , n$. In each next cycle $n + k$, $D_oT$ either produces a result that is equal to $D_oT(n + k)$ or does not give a result because $D_oT$ is not defined for $n + k$. In the second case, we put $D_oT(n + k) = D_oT(n + k - 1)$. Inductive reasoning shows that $D_oT$ determines a total function.

Proposition 2.1 is proved.

Another way to prove Proposition 2.1 for a finite set $X$ is to use the property that any finite language (set of words) is regular, i.e., $X$ is the language of some finite automata (Hopcroft *et al*, 2001).

**Example 2.4.** In the theory of inductive Turing machines, sets that are enumerable by inductive Turing machines are called inductively enumerable sets. Sets that are computable by inductive Turing machines are called inductively computable sets. They play an important role because they constitute the largest class of sets that are constructible/computable by super-recursive algorithms that give results in finite time of computation. This is based on the following conversion of Lemma 2.1 for inductive Turing machines.

**Proposition 2.2.** Any inductively computable set $X$ is in inductively enumerable.

Proof. We may assume that $X$ is an infinite set because for finite sets the statement of Proposition 2.2 follows form the Proposition 2.1. Let us suppose that $X$ is an inductively computable infinite set. It means that there is an inductive Turing machine $M$ that all outputs of $M$ constitute $X$. If $M$ is defined for all inputs, then $X$ is recursively enumerable and everything is proved. Otherwise, $M$ gives no output for some inputs, while for enumeration it is necessary to transform each input into some output. To remedy this deficiency and to preserve at the same time the initial output domain (the range of the function $f_M$), we build an inductive Turing machine $H$ that realizes the process that is called the *bidiagonal covering* of $M$. Informally, the *bidiagonal covering* means that $H$ imitates all acts of computation that $M$ performs with all possible inputs. In particular, $H$ outputs all results of $M$ and only these results. The schema of the inductive Turing machine $H$ is similar to the schema of the Turing machine $D$, which is given in Figure 1. The main difference is that $H$ has a countable number of output tapes and writes each potential output of $M$ on a separate tape. Some of these tapes give no result, while others generate all outputs of $M$.

Using the *bidiagonal covering*, which is realized by the inductive Turing machine $H$, we can correspond to $M$ an inductive Turing machine $HM$ with the following property: if $M$ is defined for a number $m$, then $HM(n) = M(m)$ and if we take all smaller numbers 1, 2, … , $m-1$, then the inductive Turing machine $M$ is defined exactly for $n-1$ of them. This machine $HM$ works in the following manner. Given a number $n$ as an input, $DT$ uses $D$ to find $n$ first results of $M$. The machine $DT$ gives the last of these results as it output for the input $n$. By induction, we demonstrate that the Turing machine $DT$ is defined for any number $k < n$. These properties of the Turing machine $DT$ show that $DT$ defines a total function $f_{DT}$ and has the same range as $T$, By Definition 2.4, this means that $DT$ enumerates the infinite set $X$.

**Definition 2.5 (Acceptation).** A set $X$ is called acceptable in **A** if all elements from $X$ and only these elements are acceptable by some algorithm $A$ from **A**.

Acceptation as a mode of computer functioning may look artificial because any computer produces some output. Computer are built to give solutions to diverse problems and they have to communicate their results to users. However, acceptation is a necessary function which is implicitly or explicitly included in any computational process. Alan Turing, making analysis of computation in his breakthrough paper (1936), wrote that to do anything with an input symbol, the device has to recognize this symbol. The first step of such recognition is acceptance of those symbols with which this device works and rejection all other symbols. For example, a device that works with the binary alphabet {0, 1} will accept the symbol "1" reject the symbol "2." This is an implicit form of acceptance. A lot of engineering efforts are aimed at the development of reliable and flexible acceptance.

The explicit form of acceptance has become important when computers acquired communication functions. For example, to be safe, an e-mail system has to accept only such messages that do not contain viruses. Thus, acceptance mode becomes vital to e-mail systems.

**Example 2.5.** In the theory of finite automata, sets acceptable by finite automata are called regular languages (Hopcroft *et al*, 2001). In the theory of pushdown automata, sets acceptable by finite automata are called context-free languages (Hopcroft *et al*, 2001).

**Definition 2.6 (Weak Decision).** A set $X$ is called weakly decidable in **A** if some algorithm from **A** decides (indicates) whether an arbitrary given word/configuration $u$ belongs to the set $X$.

Informally, weak decision means separation of all elements of $X$.

**Remark 2.6.** We may assume that if an algorithm $A$ weakly decides $X$, then $A$ produces as its output 1 for each input $x$ from $X$ and only for such inputs.

**Example 2.6.** In the theory of Turing machines and recursive functions, recursively weakly decidable sets coincide with recursively enumerable sets. Later we will prove this in a more general context.

**Definition 2.7 (Codecision).** A set $X$ is called complementary decidable or, simply, codecidable in **A** if some algorithm from **A** decides (indicates) whether an arbitrary given word/configuration $u$ does not belong to the set $X$.

Informally, codecision means separation of all elements that does not belong to $X$.

**Definition 2.8 (Decision).** A set $X$ is called decidable in **A** if some algorithm from **A** decides (indicates) whether an arbitrary given word/configuration $u$ belongs to the set $X$ or not.

Informally, decision means separation of all elements into two groups: $X$ and its complement.

**Remark 2.7.** Usually, when an algorithm $A$ decides whether a given word/configuration $u$ belongs to some set $X$, it gives the result $A(u) = 1$ if $u$ belongs to $X$ and the result $A(u) = 0$ if $u$ does not belong to $X$.

There are definite dependencies between the concepts of computability, acceptability, and decidability. For example, Definitions 2.5 and 2.6 imply the following result.

**Lemma 2.2.** A set $X$ is codecidable in **A** if and only if its complement **C**$X$ is weakly decidable in **A**.

To get other dependencies, we need to define some specific algorithms and compositions of algorithms.

**Definition 2.9.** An algorithm $A$ is called rewriting if for any input $u$ it gives $u$ as the output.

**Definition 2.10.** An algorithm $A$ is called constant if for any input $u$ it either gives no result or gives one and the same the output.

Let us consider two arbitrary algorithms $A$ and $B$ and define some operations with them.

**Definition 2.11.** The *sequential composition* $A \circ B$ of the algorithm $A$ with the algorithm $B$ is an algorithm $D$ such that the result of application of $D$ to any input $u$ is equal to the result of application of $A$ to the result of application of $B$ to $u$, i.e., $D(u) = A(B(u))$.

**Definition 2.12.** The *disjunctive parallel composition* $\vee AB$ of the algorithm $A$ with the algorithm $B$ is an algorithm $D$ such that the result of application of $D$ to any input $u$ is obtained in the following way: both $A$ and $B$ are applied to $u$ at the same time and $D(u) = A(u)$ if $A$ gets its result at the same time or earlier than $B$, otherwise, $D(u) = B(u)$.

**Definition 2.13.** The *conjunctive parallel composition* $\wedge AB$ of the algorithm $A$ with the algorithm $B$ is an algorithm $D$ such that the result of application of $D$ to any input $u$ is equal to $A(u)$ and $B(u)$.

Let $P$ be a unary predicate on $N$, i.e., $P$ is a function on $N$ that takes two values 1 and 0, which have the conventional meaning: 1 means "true" and 0 means "false."

**Definition 2.14.** The *P-sequential composition* $A_P B$ of the algorithm $A$ with the algorithm $B$ is an algorithm $D$ such that the result of application of $D$ to any input $u$ is undefined when $P(B(u)) = 0$ and equal to the result of application of $A$ to the result of application of $B$ to $u$, i.e., $D(u) = A(B(u))$, when $P(B(u)) = 1$.

**Definition 2.15.** The *P-conjunctive parallel composition* $\wedge_P AB$ of the algorithm $A$ with the algorithm $B$ is an algorithm $D$ such that the result of application of $D$ to any input $u$ is equal to $A(u)$ when $P(B(u)) = 1$ and to $B(u)$ when $P(B(u)) = 0$.

**Definition 2.16.** The *P-disjunctive parallel composition* $\vee_P AB$ of the algorithm $A$ with the algorithm $B$ is an algorithm $D$ such that the result of application of $D$ to any input $u$ is obtained in the following way: both $A$ and $B$ are applied to $u$ at the same time and $D(u) = A(u)$ if $A$ gets its result at the same time or earlier than $B$ and $P(A(u)) = 1$, $D(u) = B(u)$ if $B$ gets its result at the same time or earlier than $A$ and $P(B(u)) = 1$, otherwise, $D(u)$ is undefined.

**Definition 2.17.** The *P-disjunctive sequential composition* $\vee_P AB$ of the algorithm $A$ with the algorithm $B$ is an algorithm $D$ such that the result of application of $D$ to any input $u$ is obtained in the following way: if $A$ gives a result when applied to $u$ and $P(A(u)) = 1$, then $D(u) = A(u)$; if $B$ gives a result when applied to $u$, then $D(u) = B(u)$; otherwise, $D(u)$ is undefined.

Simple schemes of automata realize these compositions. For example, the *P*-sequential composition $D = A_P B$ of the algorithm $A$ with the algorithm $B$ is presented in the Figure 2.

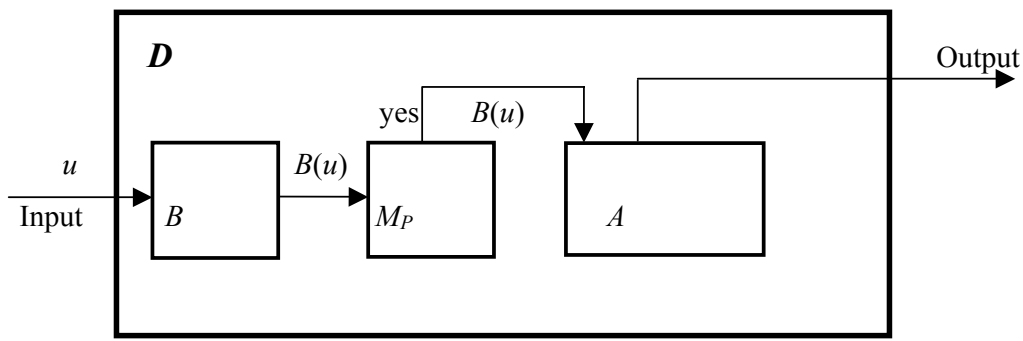

**Figure 2:** The structure of the Turing machine $D$ that realizes the $P$-sequential composition $D = A_P B$ of the algorithm $A$ with the algorithm $B$; $M_P$ is an algorithm (automaton) that realizes the predicate $P$.

Existence of such compositions in a given class of algorithms or automata results in more profound properties of this class, which represent characteristics of functioning of computers and software systems that are modeled by these algorithms or automata.

Let us assume that a class **A** of algorithms satisfies the following condition where $F$ is the predicate that takes only one value 0 for all inputs:

**(DRA)** The class **A** is closed with respect to $F$-disjunctive sequential composition of its algorithms with a rewriting algorithm.

**Remark 2.8.** When **A** contains an identity algorithm, then **A** also contains a rewriting algorithm.

**Remark 2.9.** The condition **(DRA)** is not valid for all natural classes of algorithms. The simplest examples are:

(a) a class of algorithms/automata without output;

(b) a class of algorithms/automata such that their output alphabet consists of one symbol, e.g., 1, while their output alphabet consists of more symbols, e.g., of 0, 1, 2, 3, 4, 5, 6, 7, 8, 9;

Let us assume that the class **A** satisfies one more additional condition:

**(AT)** The output alphabet of any algorithm from the class **A** contains its input alphabet.

**Proposition 2.3.** If a set $X$ is weakly decidable in **A**, then $X$ is computable in **A**.

Proof. Let $X$ be a weakly decidable set in **A**. Then by Definition 2.5, there is an algorithm $A$ from **A** decides (indicates) whether an arbitrary given word/configuration $u$ belongs to the set $X$. That is, $A(w) = 1$ when $w \in X$ and $A$ is undefined otherwise. By the initial condition, $F$-disjunctive sequential composition $D$ of $A$ with a rewriting algorithm belongs to **A**. Then by the definition of *F*-disjunctive sequential composition, *D* computes the set *X*.

Proposition 2.3 is proved.

**Remark 2.10.** In general, not every weakly decidable set is computable. For example, if we take the class of all Turing machines that give only two outputs 1 or 0, then any set that contain more than two elements is non-computable in this class, while decidable sets are also weakly decidable.

Now let us assume that a class **A** of algorithms satisfies the following condition:

**(C1)** The class **A** is closed with respect to sequential composition of its algorithms with a constant algorithm that always gives the output 1.

**Proposition 2.4.** If a set *X* is acceptable in **A**, then *X* is weakly decidable in **A**.

Indeed, let *A* belongs to **A**, and *A*(*x*) is defined for all elements *x* from *X*. Then the composition $B = A \circ C_1$ produces 1 for all inputs *x* from *X*.

Proposition 2.4 shows that acceptation is equivalent to the following two-step process: at first, make a decision whether to accept or not, then accept in the case of positive decision.

As a corollary from Propositions 2.3 and 2.4, we obtain the following result.

**Proposition 2.5.** If a set *X* is acceptable in **A**, then *X* is computable in **A**.

**Remark 2.11.** In general, not every acceptable set is computable. For example, if we take the class of all accepting finite automata, then all regular languages are acceptable by this class of algorithms, but only the empty set is computable because accepting finite automata give no output.

Let us assume that a class **A** of algorithms satisfies the following condition where *F* is the predicate that takes only one value 0 for all inputs:

**(C0)** The class **A** is closed with respect to *F*-disjunctive sequential composition of its algorithms with a constant algorithm that always gives no output.

**Proposition 2.6.** If a set *X* is weakly decidable in **A** with 1 as the indicator of membership, then *X* is acceptable in **A**.

**Proposition 2.7.** A set $X$ is weakly decidable in **A** if and only if its complement is codecidable in **A**.

In contrast to the statements of Propositions 2.3 and 2.5, this result is true by the definition of codecidability for any class of algorithms.

Let us assume that the class **A** satisfies the following condition where $F$ is the predicate that takes only one value 0 for all inputs:

**(DP)** The class **A** is closed with respect to disjunctive parallel composition of its algorithms.

**Proposition 2.8.** A set $X$ is decidable in **A** if and only if it is weakly decidable and codecidable in **A**.

Let us assume that the class **A** satisfies the following condition where $F$ is the predicate that takes only one value 0 for all inputs:

**(DCA)** The class **A** is closed with respect to $F$-disjunctive sequential composition of its algorithms with a constant algorithm.

**Proposition 2.9.** If a set $X$ is computable in **A**, then $X$ is weakly decidable in **A**.

Propositions 2.1 and 2.6 imply that computability and weak decidability are equivalent properties of sets under very weak additional conditions.

Propositions 2.3, 2.7 and 2.8 imply the following result.

**Theorem 2.1.** A set $X$ is decidable in **A** if and only if $X$ and its complement are computable in **A**.

**Remark 2.12.** In general, the statement of Theorem 2.1 may be invalid. For example, if we take the class $T_2$ of all Turing machines that give only two outputs 1 or 0, then any set that contain more than two elements is non-computable in this class, while there are many finite and infinite sets that are decidable in $T_2$.

**Proposition 2.10.** Any computable in **A** set $X$ is acceptable in **A** if the following conditions are satisfied:

**(AT)** The output alphabet of any algorithm from the class **A** is equal to its input alphabet.

**(SP)** For any algorithm $A$ from **A**, there is an algorithm $EA$ that for any input produces all outputs of $A$.

**(CA)** the class **A** realizes all comparison predicates of the form

$$\text{Comp}(w, v) = \begin{cases} 1 \text{ when } w = v; \\ 0 \text{ when } w \neq v \end{cases}$$

**(SQ)** the class **A** is closed with respect to sequential composition with comparison predicates.

**Remark 2.13.** In general, the statement of Proposition 2.10 may be invalid. For example, we can take a Turing machine $T$ that gives only for even numbers and $T(2n) = n$. Then in the class **A** that consists of $T$, only the set of all numbers is computable, while only the set of all even numbers is acceptable.

Propositions 2.5 and 2.10 imply the following result when the initial conditions of both propositions are satisfied.

**Theorem 2.2.** A set $X$ is computable in **A** if and only if it is acceptable in **A**.

It is proved in (Burgin, 1988; 1999) that there are sets such that they are inductively computable, but are not recursively computable. Thus, Theorem 2.2 implies the following result.

**Corollary 2.1.** There are sets such that they are inductively acceptable, but are not recursively acceptable.

Assuming the conditions (AT), (SP), (CA), and (SQ), we have the following theorem.

**Theorem 2.3.** A set $X$ is decidable in **A** if and only if $X$ and its complement are acceptable in **A**.

**Remark 2.14.** In general, the statement of Theorem 2.2 may be invalid. For example, we can take a Turing machine $T$ that decides some non-empty proper subset $X$ of $N$. Then in the class **A** that consists of $T$, only the set $X$ is acceptable.

It is proved in (Burgin, 1988) that there are sets such that they are inductively decidable, but are not recursively decidable. Thus, Theorems 2.1 and 2.2 imply the following result.

**Corollary 2.2.** There are sets such that they are inductively decidable, but are not recursively decidable.

## 3 Comparing Computing, Accepting, and Deciding Power of Algorithms and Programs

These concepts make possible to introduce properties of algorithmic classes that are used for their ordering. To compare different classes of algorithms by what they can do, we introduce such characteristics as computational, accepting, positive deciding, negative deciding, and deciding power of **A.**

**Definition 3.1.** The computing power of **A** is defined by the class of all sets $X$ that are computable in **A**.

**Definition 3.2.** The accepting power of **A** is defined by the class of all sets $X$ that are acceptable in **A**.

For example, the class of all Turing machines has strictly higher accepting power than the class of all pushdown automata with one stack, while the class of all pushdown automata with one stack has strictly higher accepting power than the class of all finite automata (Hopcroft *et al*, 2001).

**Definition 3.3.** The positive deciding power of **A** is defined by the class of all sets $X$ that are weakly decidable in **A**.

**Definition 3.4.** The negative deciding power of **A** is defined by the class of all sets $X$ that are codecidable in **A**.

**Definition 3.5.** The deciding power of **A** is defined by the class of all sets $X$ that are decidable in **A**.

Thus, for example, if the class of all sets *Z* that are computable in **B** includes the class of all sets *X* that are computable in **A**, then **B** has higher computing power than **A**. Comparing classes of sets in such a way, we compare power of classes of algorithms, programs, and computers.

Results from the previous section allow us to compare power of different modes of algorithms utilization and computer functioning.

Proposition 2.3 implies under the same conditions the following result.

**Theorem 3.1.** The computing power of any class **A** of algorithms is higher than the positive deciding power of **A** .

Proposition 2.4 implies under the same conditions the following result.

**Theorem 3.2.** The positive deciding power of any class **A** of algorithms is higher than the accepting power of **A** .

Proposition 2.5 implies under the same conditions the following result.

**Theorem 3.3.** The computing power of any class **A** of algorithms is higher than then accepting power of **A** .

Proposition 2.6 implies under the same conditions the following result.

**Theorem 3.4.** The accepting power of any class **A** of algorithms is higher than the positive deciding power of **A** .

As a corollary from Theorems 3.2 and 3.4, we have the following result.

**Theorem 3.5.** Accepting and positive deciding power of any class **A** of algorithms are the same if **A** satisfies the following conditions:

**(AT)** The output alphabet of any algorithm from the class **A** is equal to its input alphabet.

**(C0)** The class **A** is closed with respect to *F*-disjunctive sequential composition of its algorithms with a constant algorithm that always gives no output.

**(C1)** The class **A** is closed with respect to sequential composition of its algorithms with a constant algorithm that always gives the output 1.

As the following result demonstrates, acceptability and computability are also closely connected properties of algorithms and their classes.

**Theorem 3.6.** Computing and accepting power of any class **A** of algorithms are the same if **A** satisfies the following conditions:

**(AT)** The output alphabet of any algorithm from the class **A** is equal to its input alphabet.

**(SP)** For any algorithm $A$ from **A**, there is an algorithm $EA$ that for any input produces all outputs of $A$.

**(CA)** the class **A** realizes all comparison predicates of the form

$$\text{Comp}(w, v) = \begin{cases} 1 \text{ when } w = v; \\ 0 \text{ when } w \neq v \end{cases}$$

**(SQ)** the class **A** is closed with respect to sequential composition with comparison predicates.

**(DRA)** The class **A** is closed with respect to $F$-disjunctive sequential composition of its algorithms with a rewriting algorithm.

**(C1)** The class **A** is closed with respect to sequential composition of its algorithms with a constant algorithm that always gives the output 1.

Really, when these conditions are satisfied, any computable in **A** set $X$ is acceptable in **A** by Proposition 2.10, and by Proposition 2.5, any acceptable in **A** set $X$ is computable in **A**.

Proposition 2.8 implies under the same conditions the following result.

**Theorem 3.7.** The deciding power of any class **A** of algorithms is the intersection of its positive and negative deciding powers.

**Definition 3.6.** Two classes of algorithms are called functionally equivalent with respect to computability (acceptability, positive decidability, negative decidability or decidability) if they have the same computational (accepting, positive deciding, negative deciding or deciding, respectively) power.

For example, the class of all Turing machines is functionally equivalent to the class of all pushdown automata with two stacks (Hopcroft *et al*, 2001). For example, the class of all deterministic finite automata is functionally equivalent to the class of all nondeterministic finite automata (Hopcroft *et al*, 2001).

The recent development of the theory of algorithms implies that with respect to the way of obtaining a result, there are three modes of processing the input data (Burgin, 2000):

1. A *recursive* mode of functioning when an algorithm works with finite objects and produces its result in a finite time and stops after this, informing in such a way that the result has been obtained.

2. An *inductive* or *real-time super-recursive* mode of functioning when an algorithm works with finite objects and produces its result in a finite time, but does not inform that the result has been obtained.

3. An *unlimited super-recursive* mode of functioning when an algorithm works with infinite objects or/and produces its result in infinite time.

Turing machines and partial recursive functions work in a recursive mode, inductive inference and inductive Turing machines work in an inductive mode, and infinite time Turing machines work in an unlimited super-recursive mode.

Let $R$ be a subset of $N$, **K** be a class of Turing machines, and **H** be a class of inductive Turing machines.

**Definition 3.7.** A set $R$ is called recursive (inductive) in **K** (in **H**) if there is some (inductive) Turing machine **T** from **K** (from **H**) such that $T(r) = 1$ if $r \in R$ and $T(r) = 0$ if $r \notin R$, i.e., for any number $r$ from $N$, the machine **T** decides whether $r$ belongs to $R$ or not.

**Remark 3.1.** Recursive sets are also called decidable or recursively decidable sets and inductive sets are also called inductively decidable sets.

**Definition 3.8.** A set $R$ is called recursively (inductively) computable in **K** (in **H**) if it is computable by some recursive algorithm (some inductive Turing machine) from **K** (from **H**), i.e., some (inductive) Turing machine **T** from **K** (from **H**) computes all elements from $R$ and only such elements.

**Remark 3.2.** Recursively computable sets are also called recursively enumerable sets and inductively computable sets are also called inductively enumerable sets.

Let a class **K** (**H**) of (inductive) Turing machines is closed under the input translation, that is, for any (inductive) Turing machine **T** from **K** there is some (inductive) Turing machine **M** from **K**, such that if **M**($x$) is defined, then **M**($x$) = $x$, and **T**($x$) is defined if and only if **M**($x$) is defined.

Classes of **T** of all deterministic Turing machines, **IT** all inductive Turing machines, **NT** of all non-deterministic Turing machines, **TT** of all everywhere defined Turing machines (recursive functions), and **ITT** all everywhere defined inductive Turing machines satisfy this condition.

Then Theorem 2.1 implies the following results.

**Proposition 3.1** (Rogers, 1987). A set $R$ is recursive in **H** if and only if both $R$ and its complement are recursively computable in **H**.

The class of all Turing machines that perform computations in a polynomial time satisfies all condition from the Theorem 2.1. Thus, we can take this class as the class **H** and have the following result.

**Proposition 3.2.** A set $R$ is recursively decidable in a polynomial time if and only if both $R$ and its complement are recursively computable in a polynomial time.

The class of all Turing machines that perform computations with a polynomial space satisfies all condition from the Theorem 2.1. Thus, we can take this class as the class **H** and have the following result.

**Proposition 3.3.** A set $R$ is recursively decidable with a polynomial space if and only if both $R$ and its complement are recursively computable with a polynomial space.

**Proposition 3.4.** A set $R$ is inductive in **H** if and only if both $R$ and its complement are inductively computable in **H**.

**Remark 3.3.** As inductive inference (Blum and Blum, 1975) is realized by inductive Turing machines, Proposition 3.4 shows relations between decidability in the limit and computability in the limit.

**Remark 3.4.** Similar results may be obtained for Turing machines with advice (Karp and Lipton, 1982), Turing machines with oracles (Rogers, 1987) and persistent Turing machines (Goldin and Wegner, 1998), which model computers that are enhanced with communication functions.

## 4  Conclusion

Results comparing different modes of functioning and power of classes of algorithms, programs and computers are proved under very general axioms or conditions. This makes possible to apply these results to a vast variety of types and kinds of algorithms and their models. Such models may be structurally distinct like Turing machines and partial recursive functions. They may be defined by some restrictions inside the same class of models, e.g., polynomial time Turing machines, polynomial space Turing machines, and logarithmic time Turing machines. In its turn, comparing different models allows one to obtain relations between corresponding types of computers and software systems. For example, deterministic Turing machines model conventional computers, while nondeterministic Turing machines model quantum computers. As a result, properties of deterministic Turing machines reflect properties of conventional computers, while properties of nondeterministic Turing machines reflect properties of quantum computers.

In a similar way, recursive models of algorithms represent traditional form of computation, while super-recursive models, such as inductive Turing machines reflect pervasive computation. Results of this paper show that pervasive computation has essentially greater computing, deciding, and accepting power. As it is demonstrated in (Burgin and Shmidskii, 1996), such type of computation is rather efficient for solving many practical problems.

In addition, axiomatic approach allows one to obtain automatically many classical results of the conventional computability, which are considered in many textbooks and monographs (cf., for example, Manna, 1974; Davis and Weyuker, 1983; Hopcroft *et al*, 2001; Rogers, 1987).